\documentclass[a4paper]{amsart}

\usepackage{amsmath,amssymb,amscd,amsthm,enumerate}
\usepackage{subfigure,graphicx}

\newtheorem{Def}{Definition}[section]
\newtheorem{Th}[Def]{Theorem}
\newtheorem{Prop}[Def]{Proposition}
\newtheorem{Lem}[Def]{Lemma}
\newtheorem{Cor}[Def]{Corollary}
\newtheorem{Rem}[Def]{Remark}

\newcommand{\R}{\mathbb{R}}
\newcommand{\Z}{\mathbb{Z}}

\newcommand{\C}{\mathbb{C}}
\renewcommand{\P}{\mathbb{P}}

\newcommand{\CC}{\mathcal{C}}

\newcommand{\tp}{{\ }^{t}}
\newcommand{\al}{\alpha }
\newcommand{\be}{\beta }

\newcommand{\Ga}{\Gamma }
\newcommand{\De}{\Delta }

\newcommand{\vep}{\varepsilon }
\newcommand{\vph}{\varphi }
\newcommand{\la}{\lambda }
\newcommand{\om}{\omega }
\newcommand{\si}{\sigma }
\newcommand{\na}{\nabla }

\newcommand{\bu}{\bullet}
\newcommand{\GH}{{}_{m+1} F_m}
\newcommand{\GHE}{{}_{m+1} E_m}

\def\gh#1#2{{}_{#1}F_{#2}}
\def\para#1#2#3{\left( \begin{array}{c} #1 \\ #2 \end{array} ; #3 \right)}
\def\paraE#1#2{\left( \begin{array}{c} #1 \\ #2 \end{array} \right)}

\begin{document}

\title[Twisted (co)homology of ${}_{m+1} F_m$]
{Intersection numbers and twisted period relations 
for the generalized hypergeometric function ${}_{m+1} F_m$
}

\author[Y. Goto]
{Yoshiaki Goto}

\address{Department of Mathematics, 
Graduate School of Science,
Kobe University, 
Kobe 657-8501, Japan}
\email{y-goto@math.kobe-u.ac.jp}

\dedicatory{}

\maketitle

\begin{abstract}
We study the generalized hypergeometric function ${}_{m+1} F_m$ 
and the differential equation ${}_{m+1}E_m$ satisfied by it. 
We use the twisted (co)homology groups 
associated with an integral representation 
of Euler type.   
We evaluate the intersection numbers of some twisted cocycles 
which are defined as $m$-th exterior products of logarithmic $1$-forms. 
We also give twisted cycles corresponding to the series solutions 
to ${}_{m+1}E_m$, 
and evaluate the intersection numbers of them. 
These intersection numbers of the twisted (co)cycles lead 
twisted period relations which give 
relations for two fundamental systems of solutions to ${}_{m+1}E_m$. 
\end{abstract}

\section{Introduction}
The generalized hypergeometric function ${}_{m+1} F_m$ of
a variable $x$ with complex parameters 
$a_0 ,\ldots ,a_m ,b_1 ,\ldots ,b_m$ is defined by 
\begin{align*}
  \gh{m+1}{m} \para{a_0 ,\ldots ,a_m}{b_1 ,\ldots ,b_m}{x}
  =\sum_{n=0} ^{\infty } 
  \frac{(a_0,n) \cdots (a_m ,n)}
  {(b_1 ,n)\cdots (b_m ,n) n!} x^n, 
\end{align*}
where $b_1 ,\ldots ,b_m \not\in \{ 0,-1,-2,\ldots \}$ 
and $(c,n)=\Gamma (c+n)/\Gamma (c)$. 
This series converges in the unit disk $|x|<1$, and 
satisfies the generalized hypergeometric differential equation 
\begin{align*}
\GHE =
\GHE \paraE{a_0,\ldots ,a_m}{b_1,\ldots ,b_m}:
\left[ 
\theta \prod_{i=1}^m (\theta +b_i-1) 
-x\prod_{j=0}^m (\theta -a_j) 
\right] f(x)=0, 
\end{align*}
where $\theta =x \frac{d}{dx}$. 
The linear differential equation $\GHE$ is of rank $m+1$ 
with regular singular points $x=0,1$, and $\infty$. 
If $b_i-b_j \not\in \Z \ (0\leq i<j \leq m)$, 
a fundamental system of solutions to $\GHE$ around $x=0$ 
is given by the following $m+1$ functions: 
\begin{align}
  \begin{array}{l}
    f_0:=\GH \para{a_1 ,\ldots ,a_{m+1}}{b_1 ,\ldots ,b_m}{x} , \\
    f_{r} :=x^{1-b_r} \cdot 
    \GH \para{a_0 -b_r +1 ,\ldots ,a_m-b_r +1}
    {b_1-b_r+1 ,\ldots ,2-b_r ,\ldots ,b_m -b_r+1}{x},
  \end{array}\label{series-sol}
\end{align}
where $1\leq r \leq m$.
It is known that 
$\GH$ admits the integral representation of Euler type: 
\begin{align}
  & \GH \para{a_0,\ldots ,a_m}{b_1 ,\ldots ,b_m}{x} \label{integral} \\
  & =\prod_{i=1}^{m} \frac{\Ga(b_i)}{\Ga(a_i)\Ga(b_i-a_i)} 
  \int _{D } 
  \prod_{j=1}^{m-1} \left( t_j ^{a_j -b_{j+1}} (t_j-t_{j+1})^{b_{j+1} -a_{j+1} -1} \right) \nonumber \\
  & \hspace{40mm} \cdot t_m^{a_m -1} (1-t_1)^{b_1-a_1-1} ( 1-x t_m ) ^{-a_0} dt_1 \wedge \cdots \wedge dt_m,
  \nonumber 
\end{align}
where 
$D :=\{ (t_1 ,\ldots ,t_m) \in \R^m \mid 0<t_m <t_{m-1} <\cdots <t_1<1 \}$. 
The branch of the integrand is defined 
by the principal value for $x$ near to $0$. 

In this paper, we consider the twisted (co)homology groups associated with 
the integral representation (\ref{integral}). 
Note that the singular locus of the integrand of (\ref{integral}) is not normally crossing. 
In such a case, as is studied in \cite{KY}, the resolution of singularities 
is an effective way for the study of intersections of the twisted (co)homology groups. 
However, our singularities are so complicated for a general $m$ that 
it seems difficult to resolute them. 
To conquer this difficulty, we find a systematic method completing 
the resolution of the singular locus. 
We blow up the singular locus step by step and 
use the combinatorial structure of divisors which should exist in the complete resolution. 
The resolution of the singular locus enables us to evaluate intersection numbers
for twisted cocycles. 
We give formulas for the intersection numbers of $m$-th exterior products of 
logarithmic $1$-forms, which span the twisted cohomology group. 
For the study of the twisted homology group, we avoid the complexity of the resolution. 
We construct twisted cycles corresponding to 
the $m+1$ solutions (\ref{series-sol}) to $\GHE$ 
by using the method given in \cite{G-FC} and \cite{G-FA}. 
They are made by the bounded chambers, and their boundaries are 
canceled by different ways from the usual regularization. 
It is an advantage of our construction that we can evaluate their intersection numbers 
by the formula in \cite{Y} for a normally crossing singular locus. 
Intersection numbers of twisted homology and cohomology groups 
imply twisted period relations for two fundamental systems 
of solutions to $\GHE$ with different parameters. 
These relations are transformed into 
quadratic relations among hypergeometric series $\GH$'s. 
Since our intersection matrices are diagonal, 
it is easy to reduce the twisted period relations to quadratic relations 
among $\GH$'s. 

In \cite{Mimachi}, twisted cycles corresponding to the solutions (\ref{series-sol}) 
to $\GHE$ are obtained from real (non-bounded) chambers, 
and their intersection numbers are evaluated by the method in \cite{KY}. 
Since these cycles are scalar multiples of ours as elements of 
the twisted homology group, we give their explicit correspondence in 
Remark \ref{cf-Mimachi}. 
Twisted period relations for $\GH$ are given in \cite{Ohara} 
by the study of the intersection forms of (co)homology groups 
with coefficients in the local system of rank $m$ 
given as the solution space to $\GHE$. 
Another integral representation of $\GH$ and 
its inductive structure are used in \cite{Ohara}.

As is in \cite{BH}, the irreducibility condition of 
the differential equation $\GHE$ is 
known to be 
$a_i -b_j \not\in \Z \ (0\leq i,j \leq m)$, 
where we put $b_0:=0$ (though $b_0$ is usually defined by 1, 
we use this setting for our convenience). 
Since we use the fundamental system (\ref{series-sol}) of solutions to $\GHE$, 
we assume throughout this paper that the parameters $a_i,\ b_j$ 
satisfy the condition 
\begin{align}\label{condition}
  a_i -b_j \not\in \Z  \ (0\leq i,j \leq m),\quad 
  b_i -b_j \not\in \Z  \ (0\leq i<j \leq m). 
\end{align}

\section{Twisted (co)homology groups 
associated with the integral representation (\ref{integral})}\label{section-THG}
For twisted homology groups, twisted cohomology groups, 
and the intersection forms, 
refer to \cite{AK}, \cite{Y}, or \cite{G-FC}. 
We use the same notations as in \cite{G-FC} and \cite{G-FA}. 

In this paper, we mainly consider 
the twisted (co)homology group in \cite{G-FC} for 
$$
M:=\C ^m -\left( \bigcup_{j=1}^m (t_j=0) \cup 
\bigcup_{j=2}^m (t_{j-1}-t_j=0) \cup (1-t_1=0) \cup (1-xt_m=0)  \right) 
$$
and the multi-valued function 
\begin{align*}
  u :=& \prod_{j=1}^m t_j^{a_j-b_{j+1}} \cdot 
  \prod_{j=2}^m (t_{j-1}-t_j)^{b_j-a_j} \cdot (1-t_1)^{b_1-a_1} \cdot (1-xt_m)^{-a_0} .
\end{align*}
We put $\om :=d \log u$, where $d$ is the exterior derivative 
with respect to the variables $t_1 ,\ldots ,t_m$ 
(not to $x$ regarded as a parameter). 
The twisted cohomology group, 
that with compact support, and the twisted homology group 
are denoted by $H^k (M,\na_{\om})$, $H_c ^k (M,\na_{\om})$, and 
$H_k(\CC_{\bu} (M,u))$, respectively. 
Here, $\na_{\om}$ is the covariant differential operator defined as 
$\na_{\om} :=d+\om \wedge$. 
The expression (\ref{integral}) means that the integral 
$$
\int_{D \otimes u} u \vph_0 ,\ \ 
\vph_0  :=\frac{dt_1 \wedge \cdots \wedge dt_m }
{t_m (1-t_1) (t_1-t_2) \cdots (t_{m-1}-t_m)} 
$$
represents $\GH$ modulo Gamma factors. 
By \cite{AK} and \cite{Cho}, we have $H^k(M,\na_{\om})=0$ $(k \neq m)$, 
$\dim H^m (M,\na_{\om})=m+1$, 
and there is a canonical isomorphism 
$$
\jmath : H^m (M,\na_{\om}) \to H_c ^m (M,\na_{\om}) . 
$$
By the Poincar\'{e} duality, we have 
\begin{align*}
  & \dim H_k (\CC_{\bu} (M,u))=\dim H^k(M,\na_{\om})=0 \ \  (k \neq m),\\
  & \dim H_m (\CC_{\bu} (M,u))=\dim H^m(M,\na_{\om})=m+1. 
\end{align*}

The intersection form $I_h$ on the twisted homology groups is 
the pairing between $H_m (\CC_{\bu} (M,u))$ and $H_m (\CC_{\bu} (M,u^{-1}))$. 
The intersection form $I_c$ on the twisted cohomology groups is 
the pairing between $H_c^m (M,\na_{\om})$ and $H^m (M,\na_{-\om})$. 
By using $\jmath$, we can regard the intersection form $I_c$ 
as the pairing between $H^m (M,\na_{\om})$ and $H^m (M,\na_{-\om})$, i.e., 
$$
I_c (\psi ,\psi'):=\int_M \jmath (\psi) \wedge \psi' ,\quad 
\psi \in H^m(M,\na_{\om}) ,\ \psi' \in H^m (M,\na_{-\om}). 
$$

\section{Twisted cohomology groups and intersection numbers}
In this section, we give two systems of twisted cocycles, 
and evaluate their intersection numbers. 

We embed $M$ into the projective space $\P^m$, that is, 
we regard $M$ as the open subset of $\P^m$: 
\begin{align*}
  M=\P^m -\left( \bigcup _{j=0}^m L_j \cup \bigcup _{j=0}^m H_j  \right)
  \subset \C^m \subset \P^m ,
\end{align*}
where
\begin{align*}
  & L_j :=(T_j=0) \ (0 \leq j \leq m), \\ 
  & H_j :=(T_{j-1}-T_j =0) \ (1 \leq j \leq m) ,\ H_0 :=(T_0-xT_m =0). 
\end{align*}
By the homogeneous coordinates $T_0 ,\ldots ,T_m$, 
the multi-valued function $u$ is expressed as 
\begin{align*}
  u=T_0^{\la_0} (T_0 -xT_m)^{\mu_0} \cdot 
  \prod_{j=1}^m T_j^{\la_j} (T_{j-1} -T_j)^{\mu_j} ,
\end{align*}
where 
\begin{align*}
  & \la _j :=a_j -b_{j+1} \ (1 \leq j \leq m-1) ,\ \la_m :=a_m ,\\
  & \mu _j :=b_j -a_j \ (1 \leq j \leq m) ,\ \mu_0 :=-a_0 ,\\
  & \la_0 :=-\left( \sum_{j=1}^{m} \la_j +\sum_{j=0}^{m} \mu_j \right) =a_0 -b_1.
\end{align*}
Note that $L_0=(T_0 =0)$ is the hyperplane at infinity, 
i.e., $M \subset \C^m =\P^m -L_0 \subset \P^m$ and 
the coordinates $t_1 ,\ldots ,t_m$ on $\C^m$ are defined as $t_j =T_j /T_0$. 
Hereafter, we regard subscripts as elements in $\Z /(m+1)\Z$. 
For example, we have $a_{m+1}=a_0$, $b_{m+1}=b_0=0$, and 
$$
\la_j =a_j-b_{j+1} ,\ \mu_j =b_j -a_j \quad  (0 \leq j \leq m).
$$

Let $\ell_k$ and $h_k\ (0\leq k \leq m)$ be 
the defining linear forms of $L_k$ and $H_k$, respectively. 
We define an $m$-form on $M$ by 
$$
\phi (f_0,\ldots ,f_m):=
d\log \left( \frac{f_0}{f_1} \right) \wedge 
d\log \left( \frac{f_1}{f_2} \right) \wedge \cdots \wedge 
d\log \left( \frac{f_{m-1}}{f_m} \right) 
$$
for $f_0,\ldots ,f_m \in \{ \ell_0 ,\ldots ,\ell_m ,h_0,\ldots h_m \}$. 
We consider two systems $\{ \vph_k \}_{k=0}^m$ and $\{ \psi_k \}_{k=0}^m$ 
given as 
\begin{align*}
  & \vph_k :=\phi (h_0 ,\ldots ,h_{k-1}, \ell_{k-1} ,h_{k+1} ,\ldots h_m) , \\
  & \psi_k :=\phi (h_0 ,\ldots ,h_{k-1}, \ell_k ,h_{k+1} ,\ldots h_m). 
\end{align*}
Using the coordinates $t_j=T_j/T_0\ (1\leq j \leq m)$ of $\C^m =\P^m-L_0$, 
we have 
\begin{align*}
  \vph_0 &
  =\frac{dt_1 \wedge \cdots \wedge dt_m}{t_m(1-t_1)(t_1-t_2)\cdots (t_{m-1}-t_m)} ,\\
  \psi_0 &
  =\frac{dt_1 \wedge \cdots \wedge dt_m}{(1-t_1)(t_1-t_2)\cdots (t_{m-1}-t_m)} ,\\
  \vph_r &
  =\frac{x dt_1 \wedge \cdots \wedge dt_m}
   {t_{r-1}(1-xt_m)(1-t_1)(t_1-t_2)\cdots \widehat{(t_{r-1}-t_r)} \cdots (t_{m-1}-t_m)} ,\\
  \psi_r &
  =\frac{dt_1 \wedge \cdots \wedge dt_m}
   {t_r(1-xt_m)(1-t_1)(t_1-t_2)\cdots \widehat{(t_{r-1}-t_r)} \cdots (t_{m-1}-t_m)} ,
\end{align*}
where $1\leq r \leq m$. 
Note that the $m$-form $\vph_0$ coincides with that defined in Section \ref{section-THG}. 
\begin{Th}\label{C-intersection}
  \begin{align}
      I_c (\vph_i ,\vph_j)&=I_c (\psi_i ,\psi_j)=0  \quad  (i \neq j) , 
      \label{c-int-0} \\
      I_c (\vph_k ,\vph_k)&=(2\pi \sqrt{-1})^m
      \prod_{\substack{0\leq l \leq m\\l\neq k}}
      \frac{b_l-b_k}{(a_l-b_k)(b_l-a_l)} 
      \label{c-int-vph-self},\\
      I_c (\psi_k ,\psi_k)&=(2\pi \sqrt{-1})^m
      \prod_{\substack{0\leq l \leq m\\l\neq k}}
      \frac{a_l-a_k}{(b_l-a_k)(b_l-a_l)}
      \label{c-int-psi-self},\\
      I_c (\vph_i ,\psi_j)&=I_c (\psi_j ,\vph_i )
      =\vep_{ij} (2\pi \sqrt{-1})^m 
      \frac{(b_i -a_i)(b_j -a_j)}{(b_i -a_j)} \prod_{l=0}^{m} \frac{1}{b_l-a_l} 
      \label{c-int-vph-psi},
  \end{align}
  where 
  $$
  \vep_{ij}:=\left\{
    \begin{array}{cl}
      -1 & (i\neq j \ {\rm and}\ (i=0 \ {\rm or}\ j=0)), \\
      1 & ({\rm otherwise}).
    \end{array} \right.
  $$
\end{Th}
The following corollary follows from this theorem immediately.
\begin{Cor}\label{cor-cohomology}
  Under the condition (\ref{condition}),
  $\vph_0 ,\ldots ,\vph_m$ form a basis of $H^m (M,\na_{\om})$. 
\end{Cor}
\begin{proof}   
  Let $C:=\left( I_c (\vph_i,\vph_j) \right) _{i,j=0,\ldots,m}$
  be the intersection matrix. 
  Then we have 
  \begin{align*}
    \det (C)=(2\pi \sqrt{-1})^{m(m+1)} \prod_{l=0}^{m} \frac{1}{(b_l-a_l)^m} 
    \prod_{0 \leq i \neq j \leq m} \frac{b_i-b_j}{a_i-b_j} ,
  \end{align*}
  which does not vanish under the condition (\ref{condition}). 
\end{proof}

In the remainder of this section, we prove Theorem \ref{C-intersection}. 
According to \cite{M-form}, 
to evaluate intersection numbers, 
we have to blow up $\P^m$ so that the pole divisor of 
the pull back of $\om =d\log u$ is normally crossing. 
And we need informations of the $m$-forms 
around the points at which $m$ components of 
the pole divisor intersect. 

For $i \neq j,\ j+1$, 
let $L_{j,j+1,\dots ,i-1}$ be the exceptional divisor 
obtained by blowing up along 
$L_j \cap L_{j+1} \cap \cdots \cap L_{i-1} =(T_j=T_{j+1}=\cdots =T_{i-1}=0)$. 
The residue of the pull-back of $\om$ along $L_{j,j+1,\dots ,i-1}$ is 
\begin{align*}
  \la_{j,j+1,\dots ,i-1} 
  =\sum_{l=j}^{i-1} \la_l +\sum_{l=j+1}^{i-1} \mu_l 
  =\sum_{l=j}^{i-1} (a_l-b_{l+1}) +\sum_{l=j+1}^{i-1} (b_l-a_l) 
  =a_j-b_i
\end{align*}
(recall that the indices are regarded as elements in $\Z/(m+1)\Z$). 
Note that for example, $L_{12}$ is an exceptional divisor, 
however $L_1$ is not. 

First, we investigate the intersections of 
$L_1,\ldots ,L_m ,H_0 ,\ldots ,H_m$, and exceptional divisors 
obtained by blowing up along $L_j \cap L_{j+1} \cap \cdots \cap L_m \ (1\leq j \leq m-1)$, 
in $\C^m =\P^m -L_0$. 
By a straightforward calculation in $\P^m -L_0$, 
we obtain the following lemma. 
\begin{Lem}\label{blowup-Cm}
  We blow up $\P^m -L_0$ along 
  $$
  L_j \cap L_{j+1} \cap \cdots \cap L_{m} 
  \quad (1\leq j \leq m-1) . 
  $$
  In $\{ \vph_0 ,\ldots ,\vph_m ,\psi_0 ,\ldots ,\psi_m \}$, 
  only $\vph_0$ and $\psi_j$ have $L_{j,j+1,\ldots,m}$ 
  as a component of the pole divisor. 
  Further, we have 
  $$
  H_k \cap L_{j,j+1,\ldots,m} =\emptyset 
  \Longleftrightarrow k=0 \ {\rm or}\ k=j.
  $$
\end{Lem}

Second, we describe the all intersections of
$L_0,\ldots ,L_m ,H_0 ,\ldots ,H_m$, and exceptional divisors. 
We use the combinatorial structure of them, which arise from 
the similarity between the expression of $u$ on $\P^m-L_0$ 
and that on $\P^m -L_k$. 
\begin{Lem}\label{blow-up}
  After blowing up along all $L_j \cap L_{j+1} \cap \cdots \cap L_{i-1}$ 
  (with a suitable order), 
  the pole divisor of the pull-back of $\om$ is normally crossing. 
  Let $\Phi_k$ (resp. $\Psi_k$) be a set consisting of the components 
  of the pole divisor of the pull-back of 
  $\vph_k$ (resp. $\psi_k$). Then we have  
  \begin{align*}
    \Phi_k &=\{
    H_{k+1},\ H_{k+2},\ldots ,H_{k-1},\ 
    L_{k+1, k+2,\ldots ,k-1},\ L_{k+2,\ldots ,k-1}, 
    \ldots ,\ L_{k-2,k-1},\ L_{k-1} \} ,\\ 
    \Psi_k &= \{
    H_{k+1},\ H_{k+2},\ldots ,H_{k-1},\ 
    L_{k},\ L_{k,k+1}, 
    \ldots ,\ L_{k,k+1,\ldots ,k-3},\ L_{k,k+1 \ldots ,k-3, k-2}  \}.
  \end{align*}
  Moreover, we have 
  $$
  H_k \cap L_{j,j+1,\ldots,i-1} =\emptyset 
  \Longleftrightarrow k=i \ {\rm or}\ k=j.
  $$
\end{Lem}
\begin{proof}
  Recall that $u$ is expressed as 
  \begin{align*}
    u=t_1^{\la_1} \cdots t_m^{\la_m} \cdot 
    (1-xt_m)^{\mu_0} (1-t_1)^{\mu_1} 
    \cdot (t_1-t_2)^{\mu_2} \cdots (t_{m-1}-t_m)^{\mu_m} 
  \end{align*}
  on $\C^m =\P^m -L_0$ 
  (with coordinates $t_j=T_j/T_0 \ (1\leq j \leq m)$). 
  On the other hand, on $\P^m-L_k$, 
  it is expressed as 
  \begin{align*}
    u=&s_0^{\la_0} \cdots s_{k-1}^{\la_{k-1}} s_{k+1}^{\la_{k+1}} \cdots s_m^{\la_m} \\
    &\cdot (s_0-xs_m)^{\mu_0}(s_0-s_1)^{\mu_1} \cdots (s_{k-1}-1)^{\mu_{k}}(1-s_{k+1})^{\mu_{k+1}} 
    \cdots (s_{m-1}-s_m)^{\mu_m} \\
    =& s_{k+1}^{\la_{k+1}} \cdots s_m^{\la_m} s_0^{\la_0} \cdots s_{k-1}^{\la_{k-1}} 
    \cdot (s_{k-1}-1)^{\mu_k} (1-s_{k+1})^{\mu_{k+1}} \\
    &\cdot (s_{k+1}-s_{k+2})^{\mu_{k+2}} \cdots (s_{m-1}-s_m)^{\mu_m}
    (s_0-xs_m)^{\mu_0}(s_0-s_1)^{\mu_1} \cdots (s_{k-2}-s_{k-1})^{\mu_{k-1}} ,
  \end{align*}
  in terms of coordinates $s_j=T_j/T_k \ (0\leq j \leq m,\ j\neq k)$. 
  Thus,
  $$
  L_{k+1} ,\ldots ,L_m ,L_0 ,\ldots ,L_{k-1} \ {\rm and} \ 
  H_{k} ,\ldots ,H_m ,H_0 ,\ldots ,H_{k-1} 
  $$
  in $\P^m -L_k$ behave similarly to 
  $$
  L_1 ,\ldots ,L_m \ {\rm and} \ 
  H_0 ,\ldots ,H_m 
  $$
  in $\P^m -L_0$. 
  Then we obtain this lemma by Lemma \ref{blowup-Cm}. 
\end{proof}
\begin{Rem}
  The slight differences come from the signs of 
  $s_0-xs_m$ and $s_{k-1}-1$ at the intersection points. 
  As mentioned below, these differences make 
  complexity of $\vep_{ij}$. 
\end{Rem}


In particular, we have $\# \Phi_k =\# \Psi_k =2m$. 
We put 
\begin{align*}
  \Phi_k^{(m)}:=\left\{ 
    \{ D_1,\ldots,D_m \} \subset \Phi_k \mid 
    D_i \neq D_j \ (i\neq j),\ 
    D_1 \cap \cdots \cap D_m \neq \emptyset 
  \right\}
\end{align*}
($\Psi_k^{(m)}$ is also defined in a similar way). 
Then Lemma \ref{blow-up} implies 
\begin{align*}
  \Phi_k^{(m)} =&\left\{  
    \{ H_p \}_{p\in I} \cup \{ L_{q,q+1,\ldots ,k-1} \}_{q\not\in I}
    \mid I\subset \{ k+1, k+2,\ldots ,m,0, \ldots ,k-1 \}
  \right\} ,\\
  \Psi_k^{(m)} =&\left\{  
    \{ H_p \}_{p\in I} \cup \{ L_{k,k+1,\ldots ,q-1} \}_{q\not\in I}
    \mid I\subset \{ k+1, k+2,\ldots ,m,0, \ldots ,k-1 \}
  \right\} .
\end{align*}
Finally, we evaluate the intersection numbers of 
$\vph_i$'s and $\psi_j$'s 
by using results in \cite{M-form}. 
\begin{proof}[Proof of Theorem \ref{C-intersection}]
  First, we obtain (\ref{c-int-0}), since it is clear that 
  $$
  \Phi_i^{(m)} \cap \Phi_j^{(m)} 
  =\Psi_i^{(m)} \cap \Psi_j^{(m)} =\emptyset \quad (i \neq j). 
  $$
  Second, we have 
  \begin{align*}
    I_c (\vph_k,\vph_k) &=(2\pi \sqrt{-1})^m 
    \sum_{I\subset \{ k+1, k+2, \ldots ,k-1 \}}
    \prod_{i\in I} \frac{1}{\mu_i} 
    \cdot \prod_{j\not\in I} \frac{1}{\la_{j,j+1,\ldots ,k-1}} \\
    &=(2\pi \sqrt{-1})^m 
    \sum_{I\subset \{ k+1, k+2, \ldots ,k-1 \}}
    \prod_{i\in I} \frac{1}{b_i-a_i} 
    \cdot \prod_{j\not\in I} \frac{1}{a_j-b_k}.
  \end{align*}
  By induction on $m$, 
  we can show that 
  \begin{align*}
    \sum_{I\subset \{ k+1, k+2, \ldots ,k-1 \}}
    \prod_{i\in I} \frac{1}{b_i-a_i} 
    \prod_{j\not\in I} \frac{1}{a_j-b_k}
    =\prod_{\substack{0\leq l \leq m \\ l\neq k}}
    \frac{b_l-b_k}{(a_l-b_k)(b_l-a_l)}, 
  \end{align*}
  which implies (\ref{c-int-vph-self}). 
  The equality (\ref{c-int-psi-self}) can be shown 
  in a similar way. 
  Finally, we prove (\ref{c-int-vph-psi}). 
  Because of 
  \begin{align*}
    \Phi_i^{(m)} \cap \Psi_i^{(m)} 
    &= \left\{ \{ H_0 ,\ldots ,H_m \} -\{ H_i \} \right\}, \\
    \Phi_i^{(m)} \cap \Psi_j^{(m)} 
    &= \left\{ \{ H_0 ,\ldots ,H_m ,L_{j,j+1 ,\ldots ,i-1} \} 
      -\{ H_i,H_j \} \right\} \ (i\neq j), 
  \end{align*}
  we have 
  \begin{align*}
    I_c (\vph_i ,\psi_i)& 
    =\vep'_{ii} \cdot (2\pi \sqrt{-1})^m 
    \prod_{l\neq i} \frac{1}{\mu_l} 
    =\vep'_{ii} \cdot (2\pi \sqrt{-1})^m 
    \prod_{l\neq i} \frac{1}{b_l-a_l}, \\
    I_c (\vph_i ,\psi_j)& 
    =\vep'_{ij} \cdot (2\pi \sqrt{-1})^m \cdot 
    \frac{1}{\la_{j,j+1,\dots ,i-1}}
    \prod_{l\neq i,j} \frac{1}{\mu_l} \\
    &=\vep'_{ij} \cdot (2\pi \sqrt{-1})^m \cdot 
    \frac{1}{a_j -b_i}
    \prod_{l\neq i,j} \frac{1}{b_l-a_l} \quad (i \neq j) ,
  \end{align*}
  where $\vep'_{ij} =\pm 1$. Let us show that
  \begin{align}  
    \vep'_{ij}=\left\{
      \begin{array}{cl}
        1 & (i=0 \ {\rm or} \ j=0 \ {\rm or} \ i=j), \\
        -1& ({\rm otherwise}).
      \end{array} \right.
    \label{proof-vep}
  \end{align}
  When we evaluate the intersection number $I_c (\vph_i ,\psi_j)$, 
  it is sufficient to consider blowing up along only 
  $$
  L_{i+1}\cap L_{i+2}\cap \cdots \cap L_{i-1}, \ 
  L_{i+2}\cap \cdots \cap L_{i-1}, \ldots ,\ 
  L_{i-2}\cap L_{i-1} 
  $$
  in the coordinate system of $\P^m-L_i$, since 
  the pole divisor of $\vph_i$ is normally crossing 
  after this blowing-up process. 
  Put $\Phi_i^{(m)} \cap \Psi_j^{(m)}=\{\{ G_1,\ldots ,G_m \}\}$, 
  and let $g_l$ be the defining linear forms of $G_l$. 
  By taking appropriate coordinates $t'_1 ,\ldots ,t'_m$, 
  we express $\vph_i \cdot \prod_l g_l$ 
  and $\psi_j \cdot \prod_l g_l$ around 
  the intersection point $G_1 \cap \cdots \cap G_m$ explicitly. 
  \begin{enumerate}[(i)]
  \item For $i=j=0$, we have 
    \begin{align*}
      & g_l=1-t'_l \ (1\leq l \leq m), \\ 
      & \vph_0 \cdot \prod_l g_l =\frac{dt'_1 \wedge \cdots \wedge dt'_m}
      {t'_1 \cdots t'_m},\quad 
      \psi_0 \cdot \prod_l g_l =dt'_1 \wedge \cdots \wedge dt'_m ,
    \end{align*}
    and the intersection point $G_1 \cap \cdots \cap G_m$ is 
    expressed as 
    \begin{align*}
      t'_l =1 \ (1\leq l \leq m). 
    \end{align*}
  \item For $i=0$ and $j \neq 0$, we have 
    \begin{align*}
      & g_j=t'_j ,\ g_l=1-t'_l \ (l\neq j), \\
      & \vph_0 \cdot \prod_l g_l =\frac{dt'_1 \wedge \cdots \wedge dt'_m}
      {t'_1 \cdots \widehat{t'_j} \cdots t'_m (1-t'_j)},\quad 
      \psi_j \cdot \prod_l g_l =\frac{dt'_1 \wedge \cdots \wedge dt'_m}
      {1-xt'_1 \cdots t'_m },
    \end{align*}
    and the intersection point $G_1 \cap \cdots \cap G_m$ is 
    expressed as  
    \begin{align*}
      t'_j =0 ,\quad t'_l =1 \ (l \neq j). 
    \end{align*}
  \item For $i\neq 0$ and $j=i$, we have 
    \begin{align*}
      & g_{m+1-i}=t'_{m+1-i} -x,\ 
      g_l=1-t'_l \ (l \neq m+1-i ), \\
      & \vph_i \cdot \prod_l g_l =x\cdot \frac{dt'_1 \wedge \cdots \wedge dt'_m}
      {t'_1 \cdots t'_m},\quad 
      \psi_0 \cdot \prod_l g_l =dt'_1 \wedge \cdots \wedge dt'_m ,
    \end{align*}
    and the intersection point $G_1 \cap \cdots \cap G_m$ is 
    expressed as 
    \begin{align*}
      t'_{m+1-i} =x ,\quad t'_l =1 \ (l \neq m+1-i). 
    \end{align*}
  \item For $i\neq 0$ and $j=0$, we have 
    \begin{align*}
      & g_{m+1-i}=t'_{m+1-i},\ g_l=1-t'_l \ (l \neq  m-i+1), \\
      & \vph_i \cdot \prod_l g_l =x\cdot \frac{dt'_1 \wedge \cdots \wedge dt'_m}
      {t'_1 \cdots \widehat{t'_{m+1-i}} \cdots t'_m (t'_{m+1-i}-x)},\\ 
      & \psi_0 \cdot \prod_l g_l =\frac{dt'_1 \wedge \cdots \wedge dt'_m}
      {t'_1 \cdots t'_m-1}
    \end{align*}
    and the intersection point $G_1 \cap \cdots \cap G_m$ is 
    expressed as 
    \begin{align*}
      t'_{m+1-i} =0 ,\quad t'_l =1 \ (l \neq m+1-i). 
    \end{align*}
  \item For $i\neq 0$ and $j \neq 0,i$, we have 
    \begin{align*}
      & g_{j-i} =t'_{j-i},\ g_{m+1-i}=t'_{m+1-i}-x ,\ 
      g_l =1-t'_l \ (l\neq j-i,m+1-i),  \\
      & \vph_i \cdot \prod_l g_l =x\cdot \frac{dt'_1 \wedge \cdots \wedge dt'_m}
      {t'_1 \cdots \widehat{t'_{j-i}} \cdots t'_m (1-t'_{j-i})},\quad  
      \psi_0 \cdot \prod_l g_l =\frac{dt'_1 \wedge \cdots \wedge dt'_m}
      {t'_1 \cdots t'_m-1}
    \end{align*}
    (note that if $j<i$ then we regard $j-i$ as $m+1+j-i$), 
    and the intersection point $G_1 \cap \cdots \cap G_m$ is 
    expressed as  
    \begin{align*}
      t'_{j-i} =0,\quad t'_{m+1-i} =x ,\quad t'_l =1 \ (l \neq j-i, m+1-i). 
    \end{align*}
  \end{enumerate}
  Hence we have (\ref{proof-vep}), and 
  complete the proof of (\ref{c-int-vph-psi}).
\end{proof}

\section{Twisted homology groups and intersection numbers}\label{section-cycle}
In this section, we construct $m+1$ twisted cycles in $M$ corresponding to 
the solutions (\ref{series-sol}) to $\GHE$.

For $0 \leq k \leq m$, we set 
\begin{align*}
  M _k :=\C ^m 
  -\left( \bigcup _{j=1}^m \Bigl( (z_j =0) \cup (1-z_j=0) \Bigr) 
    \cup \Bigl( z_k-x\prod_{j\neq k}z_j =0 \Bigr) \right) ,
\end{align*}
where $z_1,\ldots, z_m$ are coordinates of $\C^m$. 
Let $u_k$ and $\phi_k$ be 
a multi-valued function and an $m$-form on $M_k$ defined as 
\begin{align*}
  & u_k :=\prod _{j \neq k}  z_j ^{a_j -b_k} (1-z_j)^{b_j-a_j} 
  \cdot z_k^{a_k} (1-z_k)^{-a_0} \left( z_k-x\prod_{j\neq k} z_j \right) ^{b_k -a_k}, \\ 
  & \phi_k  :=\frac{dz_1 \wedge \cdots \wedge dz_m}
  {z_k \cdot \prod _{j \neq k} (1-z_j) \cdot (z_k -x \prod_{j\neq k} z_j)} ,
\end{align*}
respectively. 
Here, we regard $z_0$ as $1$; we have 
\begin{align*}
  & z_0-x\prod_{j\neq k}z_j =1-x\prod_{i=1}^m z_i , \\
  & u_0=\prod_{i=1}^m z_i ^{a_i} (1-z_i)^{b_i-a_i} 
  \cdot \left( 1-x \prod_{i=1}^m z_i  \right) ^{-a_0} ,\quad
  \phi_0 =\frac{dz_1 \wedge \cdots \wedge dz_m}{\prod_{i=1}^m \bigl( z_i (1-z_i) \bigr)} .
\end{align*}
We construct a twisted cycle $\tilde{\De }_k$ loaded by $u_k$ in $M_k$. 
Let $x$ and $\vep$ be positive real numbers satisfying 
$$
\vep <\frac{1}{2},\quad 
x <\frac{\vep}{(1+\vep)^{m-1}}  
$$
(for example, if 
$$
\vep =\frac{1}{3} ,\ 0<x<\frac{1}{3} \cdot \left( \frac{3}{4} \right)^{m-1} ,
$$
this condition holds). 
Thus the direct product 
\begin{align*}
  \sigma_k :=\left\{ (z_1 ,\ldots ,z_m )\in \R^m \mid 
    \vep \leq z_k \leq 1-\vep \ (1\leq r \leq m)
  \right\} 
\end{align*}
of $m$ intervals 
is contained in the bounded domain 
$$
\left\{ (z_1 ,\ldots ,z_m )\in \R^m \ \Big| \  0<z_j<1, \ 
z_k > x \cdot \prod_{j \neq k} z_j \right\} \subset (0,1)^m.
$$
The orientation of $\sigma _k$ 
is induced from the natural embedding $\R^m \subset \C^m$. 

By using the $\vep$-neighborhoods of 
$C_1 :=(z_1 =0)$, $\ldots$, $C_m :=(z_m =0)$, $C_{m+1} :=(1-z_1=0)$, $\ldots$, $C_{2m} :=(1- z_m =0)$, 
we construct a twisted cycle $\tilde{\De}_k$ 
from $\sigma _k \otimes u_k$ in a similar way in \cite{G-FC} and \cite{G-FA}. 
If $k\neq 0$, 
we have to consider the difference of branches of
$$
z_k^{a_k}  \left( z_k-x\prod_{j\neq k} z_j \right) ^{b_k -a_k}
$$
at the ending and starting points of a circle surrounding $C_k$. 
Indeed, for fixed positive real numbers $z_j \ (j\neq k)$, 
the solution $x\prod_{j\neq k} z_j$ of the equation 
$z_k-x\prod_{j\neq k} z_j =0$ belongs to $\R$ and satisfies 
\begin{align*}
  x \cdot \prod_{j \neq k} z_j < x <\vep .
\end{align*}
Thus, the difference is 
$\exp (2\pi \sqrt{-1}a_k) \exp (2\pi \sqrt{-1}(b_k -a_k))$ and 
the exponent about this contribution is 
$$
a_k +(b_k -a_k) =b_k.
$$
The exponents about the contributions 
of the circles surrounding $C_{m+k}$, $C_j$, $C_{m+j}$ $(j \neq k)$ 
are simply 
$$
-a_0 ,\ a_j -b_k ,\ b_j -a_j ,
$$ 
respectively. 

\begin{Rem}
  If $k=0$, the exponents about the contributions 
  of the circles surrounding $C_i$, $C_{m+i}$ $(1\leq i \leq m)$ 
  are simply $a_i -b_0=a_i$, $b_i -a_i$, respectively. 
  Since $\overline{(0,1)^m} \cap \{ z \mid 1-x\prod_i z_i=0 \} =\emptyset$, 
  the twisted cycle $\tilde{\De}_0$ is the usual regularization of $(0,1)^m \otimes u_0$. 
\end{Rem}

\begin{Prop}\label{series-cycle}
\begin{align*}
\int _{\tilde{\De }_0} u_0 \phi_0 
=& \prod_{i=1}^m \frac{\Ga(a_i) \Ga(b_i-a_i)}{\Ga(b_i)} 
\cdot \GH \para{a_0,\ldots ,a_m}{b_1,\ldots ,b_m}{x} ,\\
\int _{\tilde{\De }_r} u_r \phi_r
=& \frac{\Ga(b_r-1) \Ga(1-a_0)}{\Ga(b_r-a_0)} 
\cdot \prod _{\substack{1\leq j \leq m \\ j \neq r}} \frac{\Ga(a_j -b_r +1) \Ga(b_j -a_j)}{\Ga(b_j -b_r +1)} \\
& \cdot
\GH \para{a_0 -b_r +1 ,\ldots ,a_m-b_r +1}{b_1-b_r+1 ,\ldots ,2-b_r ,\ldots ,b_m -b_r+1}{x} 
\quad (1\leq r \leq m).
\end{align*}
\end{Prop}
\begin{proof}
  In a similar way to Proposition 4.3 of \cite{G-FC} or Proposition 4.3 of \cite{G-FA}, 
  we can show this proposition by 
  expanding the left hand sides with respect to $x$. 
  Note that we use the equalities 
  \begin{align*}
    & \int_{\tilde{\De }_0} \prod_{i=1}^m z_i^{a_i+n-1} (1-z_i)^{b_i-a_i-1}dz
    =\prod_{i=1}^m \frac{\Ga (a_i+n) \Ga (b_i-a_i)}{\Ga (b_i+n)},  \\
    & \int _{\tilde{\De }_r} \prod_{j \neq r} z_j ^{a_j -b_r +n} (1-z_j)^{b_j-a_j-1}  
    \cdot z_r^{b_r-2-n} (1-z_r)^{-a_0} dz \\
    &=\prod_{j \neq r} \frac{\Ga(a_j-b_r+n+1)\Ga(b_j-a_j)}{\Ga(b_j-b_r+n+1)} 
    \cdot \frac{\Ga(b_r-1-n)\Ga(1-a_0)}{\Ga(b_r-a_0-n)},
  \end{align*}
  for a natural number $n$ and $1\leq r \leq m$. 
  The second equality follows from the fact that 
  the twisted cycle $\tilde{\De }_r$ of the integral can be 
  identified with the usual regularization of the domain $(0,1)^m$ 
  loaded by the multi-valued function 
  $$
  \prod_{j \neq r} z_j ^{a_j -b_r +n} (1-z_j)^{b_j-a_j-1}  
  \cdot z_r^{b_r-2-n} (1-z_r)^{-a_0}
  $$
  on $\C^m -\bigcup_{j=1}^m \bigl( (z_j=0) \cup (1-z_j =0) \bigr)$. 
\end{proof}

We define a bijection $\iota _k :M_k \rightarrow M$ by
\begin{align*}
  & \iota _0 (z_1 ,\ldots ,z_m ):=(t_1 ,\ldots ,t_m );\ 
  t_s=\prod_{i=1}^s z_i ,\\
  & \iota _r (z_1 ,\ldots ,z_m ):=(t_1 ,\ldots ,t_m );\ 
  t_s =\prod_{j=1}^s z_j \ (s<r),\ 
  t_s =\frac{z_r}{x\cdot \prod_{j=s+1}^m z_j} \ (s\geq r) ,
\end{align*}
where $1 \leq r \leq m$. 
We also define branches of the multi-valued function $u$ 
on real chambers in $M$. 
Let $D_r \subset \R^m \ (1\leq r \leq m)$ be the chamber defined by 
\begin{align*}
  t_j>0\ (1\leq j \leq m),\ 1-xt_m >0,\ 
  t_{j-1}-t_j>0 \ (j\neq r),\ t_{r-1}-t_r<0 ,
\end{align*}
where we regard $t_0$ as $1$. 
On $D_r$, 
the arguments of the factors of $u$ are given as follows. 
\begin{eqnarray*}
  \begin{array}{|c|c|c|c|}
    \hline
    t_j & 1-xt_m & t_{j-1}-t_j \ (j\neq r) & t_{r-1}-t_r \\ \hline
    0 & 0 & 0 & -\pi \\ \hline
  \end{array}
\end{eqnarray*}
Recall that 
on $D =\{ (t_1 ,\ldots ,t_m) \in \R^m \mid 0<t_m <t_{m-1} <\cdots <t_1<1 \}$, 
all of the arguments of the factors of $u$ are $0$. 

\begin{Th}\label{solutions-cycle}
We define a twisted cycle $\De _k$ in $M$ by 
\begin{align*}
  \De _k :=(\iota _k )_{*} (\tilde{\De }_k) .
\end{align*}
Then we have
\begin{align*}
  &\int_{\De_0} u \vph_0 
  =\prod_{i=1}^m \frac{\Ga(a_i) \Ga(b_i-a_i)}{\Ga(b_i)} \cdot f_0 ,\\
  &\int_{\De_r} u \vph_0 
  =e^{-\pi \sqrt{-1} (b_r-a_r-1)} 
  \frac{\Ga(b_r-1) \Ga(1-a_0)}{\Ga(b_r-a_0)} 
  \cdot \prod _{\substack{1\leq j \leq m \\ j \neq r}} \frac{\Ga(a_j -b_r +1) \Ga(b_j -a_j)}{\Ga(b_j -b_r +1)}
  \cdot f_r ,
\end{align*}
where $1\leq r \leq m$. 
\end{Th}
\begin{proof}
By pulling back $u\vph_0$ under $\iota_0$, 
we can show the first claim. 
We prove the second one. 
On $\De_r$, we have 
\begin{align*}
  u=&e^{-\pi \sqrt{-1} (b_r-a_r)}
  (t_r-t_{r-1})^{b_r-a_r} (1-xt_m)^{-a_0} \\
  &\cdot \prod_{j=1}^{m} t_j^{a_j-b_{j+1}} \cdot 
  \prod_{\substack{1\leq j \leq m \\ j\neq r}} (t_{j-1}-t_j)^{b_j-a_j} ,
\end{align*}
where the argument of each factor is zero on $\iota_r (\si_r) \subset D_r$. 
We consider the pull back of $u \vph_0$ under $\iota_r$: 
\begin{align*}
  u(\iota_r (z))
  =&
  e^{-\pi \sqrt{-1} (b_r-a_r)} \cdot x^{-b_r} \cdot u_r (z), \\
  \iota_r^{*} \vph_0
  =&
  -x\cdot \phi_r . 
\end{align*}
By Proposition \ref{series-cycle}, we thus have 
$$
\int_{\De_r} u \vph_0 
=-e^{-\pi \sqrt{-1} (b_r-a_r)} x^{1-b_r} \int _{\tilde{\De }_r} u_r \phi_r 
=e^{-\pi \sqrt{-1} (b_r-a_r-1)}  \cdot (\Ga {\textrm -}{\rm factors})\cdot f_r .
$$ 
\end{proof}

\begin{Rem}
  For $1\leq r \leq m$, 
  the twisted cycle $\De_r$ is different from 
  the regularization of $D_ r \otimes u$ 
  as elements in $H_m (\CC_{\bu} (M,u))$. 
\end{Rem}
\begin{Rem}
  Let 
  $\iota'_r :M_r \to M_0 \ (1\leq r \leq m)$ be the map defined as 
  $$
  \iota'_r (z_1 ,\ldots ,z_m ):=(w_1 ,\ldots ,w_m );\ 
  w_r =\frac{z_r}{x\prod_{j\neq r} z_j} ,\ w_s =z_s \ (s\neq r).
  $$
  Then it is easy to see that 
  $\iota_r =\iota_0 \circ \iota'_r$. 
\end{Rem}

The replacement $u \mapsto u^{-1} =1/u$ and the construction same as $\De _k$ give 
the twisted cycle $\De_k ^{\vee}$ 
which represents an element in $H_m (\CC_{\bu} (M,u^{-1}))$. 
We obtain the intersection numbers of the twisted cycles $\{ \De _k \} _{k=0}^m$ and 
$\{ \De _k^{\vee} \}_{k=0}^m$. 

\begin{Th}\label{H-intersection}
  \begin{enumerate}[(i)]
  \item For $k\neq l$, 
    we have $I_h (\De_k ,\De_l ^{\vee}) =0$. 
  \item The self-intersection numbers of $\De_k$'s are as follows:  
    \begin{align*}
      I_h (\De _0 ,\De _0^{\vee})
      &= \prod_{i=1}^{m} \frac{\al_i (1-\be_i)}{(1-\al_i)(\al_i -\be_i)} ,\\
      I_h (\De _r ,\De _r^{\vee})
      &= \prod_{\substack{1\leq j \leq m \\ j\neq r}} 
      \frac{\al_j (\be_r-\be_j)}{(\be_r-\al_j)(\al_j -\be_j)} 
      \cdot \frac{\al_0-\be_r}{(1-\be_r)(\al_0-1)} \quad (1 \leq r \leq m),
    \end{align*}
    where $\al_j :=e^{2\pi \sqrt{-1}a_j} ,\ \be_j :=e^{2\pi \sqrt{-1}b_j}$. 
  \end{enumerate}
\end{Th}
\begin{proof}
This theorem can be also shown similarly to 
Theorem 4.6 of \cite{G-FC} or Theorem 4.6 of \cite{G-FA}. 
\end{proof}
\begin{Cor}\label{cor-homology}
  Under the condition (\ref{condition}), the twisted cycles 
  $\De_0 ,\ldots ,\De_m$ form a basis of $H_m (\CC_{\bu} (M,u))$
\end{Cor}
\begin{proof}
  The determinant of the intersection matrix 
  $H:=\left( I_h (\De_i,\De_j) \right)_{i,j=0,\ldots ,m}$ 
  does not vanish. 
\end{proof}
\begin{Rem}\label{cf-Mimachi}
  In Section 3 of \cite{Mimachi}, there are the twisted cycles 
  $D_1^{(0)} ,\ldots ,D_m^{(0)} ,D_{m+1}^{(0)}$ 
  which correspond to the solutions $f_1 ,\ldots ,f_m ,f_0$, respectively. 
  By the variable change 
  $$
  p: (t_1 ,\ldots ,t_m) \longmapsto 
  \left( \frac{1}{t_1},\ldots ,\frac{1}{t_m} \right), 
  $$
  our integral representation (\ref{integral}) coincides with 
  that in \cite{Mimachi}. 
  It is easy to see that 
  \begin{align*}
    \De_0 =(-1)^m p_*(D_{m+1}^{(0)}) ,\quad 
    \De_r =(-1)^m \frac{\be_r -\al_0}{\al_0 (\be_r -1)}
    p_*(D_r^{(0)}) \quad (1 \leq r \leq m) 
  \end{align*}
  as elements in $H_m (\CC_{\bu} (M,u))$. 
\end{Rem}

\section{Twisted period relations} \label{section-TPR}
The compatibility of the intersection forms and the pairings obtained by integrations 
(see \cite{CM}) implies twisted period relations: 
\begin{align*}
  C=\Pi_{\om} \tp H^{-1} \tp \Pi_{-\om}, 
\end{align*}
where $\Pi_{\pm \om}$ are defined as 
$$
\Pi_{\om}:=\left( \int _{\De _j} u\vph_i \right)_{i,j},\quad  
\Pi_{-\om}:=\left(\int _{\De _j ^{\vee}} u^{-1} \vph_i \right)_{i,j} , 
$$
$C$ and $H$ are the intersection matrices 
(see the proof of Corollaries \ref{cor-cohomology} and \ref{cor-homology}). 
Comparing the $(i,j)$-entries of both sides, we obtain the following theorem.  

\begin{Th} \label{TPR}
  We have 
  \begin{align}
    \label{TPR-eq}
    I_{c} (\vph_i ,\vph_j) 
    = \sum_{k=0}^m \frac{1}{I_{h} (\De _k ,\De _k^{\vee} )}
    \cdot \int _{\De _k} u\vph_i \cdot \int _{\De _k ^{\vee}} u^{-1} \vph_j .
  \end{align}
\end{Th}
By using our results, we can reduce 
the twisted period relations (\ref{TPR-eq}) 
to quadratic relations among $\GH$'s. 
We write down one of them as a corollary. 
\begin{Cor}\label{TPR-cor}
  The equality (\ref{TPR-eq}) for $i=j=0$ is reduced to 
  \begin{align*}
    \prod_{l=1}^{m} \frac{b_l-b_0}{a_l-b_0}
    =&\prod_{l=1}^{m} \frac{b_l}{a_l} 
    \cdot \GH \para{\boldsymbol{a}}{\boldsymbol{b}}{x} 
    \cdot \GH \para{-\boldsymbol{a}}{-\boldsymbol{b}}{x} \\
    &+\sum _{r=1}^{m} x^2 \cdot \frac{a_0(a_0-b_r)(b_r-a_r)}{b_r(b_r^2-1)} 
    \cdot \prod _{\substack{1\leq l \leq m \\ l\neq r}} \frac{a_l -b_r}{b_l -b_r} \\
    & \quad \quad \quad \quad \quad 
    \cdot \GH \para{\boldsymbol{a}^{r,+}}{\boldsymbol{b}^{r,+}}{x} 
    \cdot \GH \para{\boldsymbol{a}^{r,-}}{\boldsymbol{b}^{r,-}}{x}, 
  \end{align*}
  where 
  \begin{align*}
    & \boldsymbol{a} :=(a_1 ,\ldots ,a_{m+1}) ,\ 
    \boldsymbol{b} :=(b_1 ,\ldots ,b_m) , \\
    & \boldsymbol{a}^{r,\pm} :=(1,\ldots ,1)\pm (a_1-b_r,\dots ,a_{m+1}-b_r) ,\\
    & \boldsymbol{b}^{r,\pm} :=
    (1,\ldots ,1)\pm (b_1-b_r ,\dots ,\pm 1-b_r ,\dots ,b_m-b_r).
  \end{align*}
\end{Cor}
\begin{proof}
  By 
  Theorem \ref{solutions-cycle}, 
  we can express the integrals in (\ref{TPR-eq}) as 
  products of $\Ga$-factors and $\GH$. 
  By Theorems \ref{C-intersection}, \ref{H-intersection}, and 
  the formula $\Ga (w) \Ga (1-w) =\pi/\sin (\pi w)$, 
  we obtain the corollary. 
\end{proof}

\begin{Rem}
  If we assume the condition (\ref{condition}) and 
  $$
  a_i-a_j \not\in \Z \quad (0\leq i<j \leq m), 
  $$
  then $\psi_0,\ldots ,\psi_m$ also form a basis of $H^m (M,\na_{\om})$, 
  because of Theorem \ref{C-intersection}. 
  Considering $I_c (\vph_i ,\psi_j)$, $I_c (\psi_i ,\vph_j)$, or 
  $I_c (\psi_i ,\psi_j)$, we obtain other twisted period relations. 
\end{Rem}

\section*{Acknowledgments}
The author thanks Professor Keiji Matsumoto for his useful advice and constant encouragement.


\begin{thebibliography}{99}
\bibitem{AK} K. Aomoto and M. Kita, translated by K. Iohara, 
  {\it Theory of Hypergeometric Functions} 
  (Springer Verlag, New York, 2011). 

\bibitem{BH}
  F. Beukers and G. Heckman, 
  Monodromy for the hypergeometric function ${}_n F_{n-1}$, 
  {\it Invent. Math.},  {\bf 95} (1989), no. 2, 325--354. 

\bibitem{Cho} K. Cho,
  A generalization of Kita and Noumi's vanishing theorems of 
  cohomology groups of local system, 
  {\it Nagoya Math. J.}, {\bf 147} (1997), 63--69. 

\bibitem{CM}
  K. Cho and K. Matsumoto, 
  Intersection theory for twisted cohomologies and 
  twisted Riemann's period relations I,   
  {\it Nagoya Math. J.},  {\bf 139} (1995), 67--86. 

\bibitem{G-FC}
  Y. Goto, 
  Twisted cycles and
  twisted period relations
  for Lauricella's hypergeometric function $F_C$, 
  {\it Internat. J. Math.}, {\bf 24} (2013), 1350094 19pp.

\bibitem{G-FA}
  Y. Goto, 
  Twisted period relations
  for Lauricella's hypergeometric function $F_A$, preprint.

\bibitem{GM} Y. Goto and K. Matsumoto, 
  The monodromy representation and twisted period relations 
  for Appell's hypergeometric function $F_4$,
  to appear in {\it Nagoya Math. J.} 

\bibitem{KY} M. Kita and M. Yoshida, Intersection theory for twisted cycles. II. Degenerate arrangements,
  {\it Math. Nachr.}, {\bf 168} (1994), 171--190. 

\bibitem{M-form} K. Matsumoto, 
  Intersection numbers for logarithmic $k$-forms, 
  {\it Osaka J. Math.}, {\bf 35} (1998), 873--893. 

\bibitem{Mimachi} K. Mimachi,
  Intersection numbers for twisted cycles and the connection problem associated with 
  the generalized hypergeometric function ${}_{n+1}F_n$, 
  {\it Int. Math. Res. Not. IMRN}, (2011), no. 8, 1757--1781. 

\bibitem{Ohara} K. Ohara, Y. Sugiki and N. Takayama, 
  Quadratic relations for generalized hypergeometric functions ${}_p F_{p-1}$, 
  {\it Funkcial. Ekvac.}, {\bf 46} (2003), no. 2, 213--251. 
  
\bibitem{Y} M. Yoshida, 
  {\it Hypergeometric functions, my love, -Modular interpretations of configuration spaces-} 
  (Vieweg \& Sohn, Braunschweig, 1997).
\end{thebibliography}
\end{document}